\documentclass[11pt]{article}

	\usepackage{amsfonts,amsmath,amssymb,amsthm,cite}

	\def\Z{\mathbb Z}
	\def\f{f_{-n}} \def\a{a_0}
	\def\k{\kappa_0}
	\def\S{\mathbb S}
	\def\R{\mathbb R}
		
	\def\sn{\S^{n-1}} 
	  
	\def\cK{{\cal K}} 
	\def\cH{{\cal H}} 
	\def\cB{{\cal B}} 
	 
	\def\cP{{\cal P}} 
	
	\newcommand{\bd}{\partial}
	\newcommand{\inte}{\mathop{\rm int}}
	\newcommand{\reg}{\mathop{\rm reg}}
		\newcommand{\gln}{\operatorname{GL}(n)} \newcommand{\sln}{\operatorname{SL}(n)}

	\def\as{\Omega_{\psi}}

	\def\fuer{\,\,\,\hbox{ for }\,\,\,}
	
	\def\B{ B^n}

   \newcommand{\concave}{\operatorname{\rm Conc}(0,\infty)} 
    \newcommand{\convex}{\operatorname{\rm Conv}(0,\infty)}

	\newtheorem{theorem}{Theorem} 
	\newtheorem{corollary}[theorem]{Corollary}
	\newtheorem{question}{Conjecture}

	\newcommand{\eqnref}[1]{(\ref{#1})}

	\title{\sc General Affine Surface Areas}
	
        \author{Monika Ludwig\,\footnote{Research supported, in part, by NSF grant DMS-0805623}}
        \date{} 
\begin{document}
\maketitle

\begin{abstract}
Two families of general affine surface areas are introduced. {Basic} properties and affine isoperimetric inequalities for these new affine surface areas as well as for $L_{\phi}$~affine surface areas are established.

\end{abstract}


Finding the {\em right\,} notion of affine surface area  was one of the first questions asked  within affine differential geometry. At the beginning of the last century, Blaschke \cite{Blaschke} and his School  studied this question and introduced equi-affine surface area -- a notion of surface area that is equi-affine in\-va\-riant, that is, $\sln$ and translation invariant. The first fundamental result regarding equi-affine surface area was the classical affine isoperimetric inequality of differential geometry \cite{Blaschke}. Numerous important results regarding equi-affine surface area were obtained in recent years (see, for example, \cite{Andrews96,Andrews99,Sheng:Trudinger:Wang,Trudinger:Wang2000,Trudinger:Wang02,Trudinger:Wang05,Trudinger:Wang06}). 
Using valuations on convex bodies,  the author and Reitzner \cite{Ludwig:Reitzner2} were able to characterize a much richer family of affine surface areas (see Theorem \ref{class}).  Classical equi-affine and centro-affine surface area as well as all $L_p$~affine surface areas for $p>0$ belong to this family of $L_{\phi}$~affine surface areas.

The present paper has two aims. The first is to establish affine isoperimetric inequalities and basic duality relations  for all $L_{\phi}$~affine surface areas. 
The second aim is to define  new general notions of affine surface area that complement $L_{\phi}$~affine surface areas and include $L_p$~affine surface areas for $p<\!-n$ and $-n<p<0$. Let $\cK_0^n$ denote the space  of convex bodies, that is, compact convex sets, in $\R^n$ that contain the origin in their interiors. Whereas $L_{\phi}$~affine surface areas are always finite and are upper semi\-continuous functionals on $\cK_0^n$, the affine surface areas of the new families are infinite for certain convex bodies including polytopes and are lower semicontinuous functionals on $\cK^n_0$. 
Basic properties and affine isoperimetric inequalities for these new affine surface areas are established. In Section \ref{conj}, it is conjectured that together with $L_\phi$ affine surface areas, these new affine surface areas constitute -- in a certain sense -- {\em all\,} affine surface areas.

For a smooth convex body $K\subset \R^n$,  equi-affine surface area is defined by
\begin{equation}\label{asadef}
\Omega(K) = \int_{\bd K} \k(K,x)^\frac{1}{n+1}\,d\mu_K(x).
\end{equation}
Here $d\mu_K(x)=  x\cdot u(K,x)\,d\cH(x)$
is the cone measure on $\bd K$,
$x \cdot u$ is the standard inner product of $x,u\in\R^n$,
$u(K,x)$ is the exterior unit normal vector to $K$ at $x\in\bd K$,
$\cH$ is the $(n-1)$-dimensional Hausdorff measure,
\begin{equation*}
\k(K,x)=\frac{\kappa(K,x)}{(x\cdot u(K,x))^{n+1}},
\end{equation*}
and $\k(K,x)$ is the Gaussian curvature of $K$ at $x$.
Note that $\k(K,x)$ is (up to a constant) just a power of the volume of the origin-centered ellipsoid osculating $K$ at $x$ and thus is an $\sln$ covariant notion. Also $\mu_K$ is an $\sln$ covariant notion. Thus $\Omega$ is easily seen to be $\sln$ invariant and it is also easily seen to be  translation invariant.  The notion of equi-affine surface area is fundamental in affine differential and convex geometry. Since many basic problems in discrete and stochastic geo\-metry are equi-affine invariant, equi-affine surface area has found numerous applications in these fields  (see, for example,  \cite{Barany97,Barany99,Gruber93b,Reitzner05}).

The extension of the definition of equi-affine surface area to general convex bodies was obtained much more recently in a series of papers \cite{Leichtweiss88,Lutwak91,Schuett:Werner90}. Since $\k(K, \cdot)$ exists $\mu_K$ a.e.\ on $\bd K$ by Aleksandrov's differentiability theorem, definition \eqnref{asadef} still can be used.  The long conjectured upper semicontinuity of equi-affine surface area (for smooth surfaces as well as for general convex surfaces)  was proved  by Lutwak \cite{Lutwak91} in 1991, that is,
$$\limsup_{j\to\infty}\, \Omega(K_j)\le \Omega(K)$$
for any sequence of convex bodies $K_j$ converging to $K$ (in the Hausdorff metric).
 Let $\cK^n$ denote the space of convex bodies in $\R^n$. Sch\"utt \cite{Schuett93} showed that  $\Omega$ is a valuation on $\cK^n$, that is,
$$\Omega(K)+\Omega(L)=\Omega(K\cup L)+\Omega(K\cap L)$$
for all $K,L\in\cK^n$ with $K\cup L\in\cK^n$. An equi-affine version of Hadwiger's celebrated classification theorem \cite{Hadwiger:V} was established in \cite{Ludwig:Reitzner}:  (up to multiplication with a positive constant) equi-affine surface area is  the unique upper semicontinuous, $\sln$ and translation invariant valuation on $\cK^n$ that vanishes on polytopes.

\goodbreak

During the past decade and a half, there has
been an explosive growth of an $L_p$ extension of the classical
Brunn Minkowski theory (see, for example, \cite{Campi:Gronchi, Chou:Wang,Ludwig:matrix,Ludwig:Minkowski, Lutwak93b, LYZ2000, LYZ2002b,LYZ2004, LZ1997,Stancu02,Stancu03}). Within this theory, $L_p$~affine
surface area is the notion corresponding to equi-affine
surface area in the classical Brunn Minkowski theory.  
Let $\cK_0^n$ denote the space of convex bodies in $\R^n$ that contain the origin in their interiors. 
For $p>1$, $L_p$~affine surface area, $\Omega_p$, was introduced by 
Lutwak \cite{Lutwak96}  and shown to be $\sln$ invariant, homogeneous of degree $q= p(n-p)/(n+p)$ (that is, $\Omega_p(t\,K)=t^q\,\Omega_p(K)$ for $t>0$), and upper semicontinuous on $\cK_0^n$.
Hug \cite{Hug96} defined $L_p$~affine surface area for every $p>0$ and obtained the following representation for $K\in\cK_0^n$:
\begin{equation}\label{lpdef}
\Omega_p(K)=\int_{\bd K} \k(K,x)^{\frac{p}{n+p}}\,d\mu_K(x).
\end{equation}
Note that $\Omega_1=\Omega$ and that $\Omega_n$ is the classical (and $\gln$ invariant) centro-affine surface area. Geometric interpretations of $L_p$~affine surface areas were obtained in \cite{Gruber93,Meyer:Werner2000, Schuett:Werner2004,Werner2007a}, and an application of $L_p$~affine surface areas to partial differential equations  is given in  \cite{Lutwak:Oliker}.

The $L_p$~affine surface areas for $p>0$ are special cases of the following family of affine surface areas introduced in \cite{Ludwig:Reitzner2}.
Let $\concave$ be the set of functions $\phi:(0,\infty)\to (0,\infty)$ such that $\phi$ is concave, $\lim_{t\to 0}\phi(t)=0$, and $\lim_{t\to\infty} \phi(t)/t=0$.  Set $\phi(0)=0$. 
For $\phi\in\concave$, we define the $L_{\phi}$~affine surface area of $K$ by
\begin{equation}\label{one}    
\Omega_{\phi}(K) = \int_{\bd K} \phi(\k(K,x))\,d\mu_K(x).
\end{equation}
The following basic properties of $L_{\phi}$~affine surface areas were established in \cite{Ludwig:Reitzner2}. Let $\cP_0^n$ denote the set of convex polytopes containing the origin in their interiors.

\begin{theorem}[\!\!\!\cite{Ludwig:Reitzner2}]\label{prop}
If  $\,\phi\in\concave$, then $\,\Omega_{\phi}(K)$ is finite for every $K\in\cK_0^n$ and
$\Omega_{\phi}(P)=0$ for every $P\in\cP_0^n$. In addition, $\Omega_{\phi}:\cK_0^n\to [0,\infty)$ is both upper semicontinuous and  an $\,\sln$ invariant valuation.
\end{theorem}
\noindent
The family of $L_{\phi}$~affine surface areas for $\phi\in\concave$ is distinguished by the following basic properties (see \cite{Ludwig:origin} and \cite{Ludwig:Reitzner2}, for characterizations of functionals that do not necessarily vanish on polytopes).
\begin{theorem}[\!\!\!\cite{Ludwig:Reitzner2}]\label{class}
If $\,\Phi:\cK_0^n\to \R$ is an upper semicontinuous and $\,\sln$ invariant valuation that vanishes on $\cP^n_0$, then there exists $\phi\in\concave$ such that
$$\Phi(K)=\Omega_{\phi}(K)$$
for every $K\in\cK_0^n$.
\end{theorem}

One of the most important inequalities of affine geometry is the classical affine isoperimetric inequality. The following theorem establishes affine isoperimetric inequalities for all $L_{\phi}$~affine surface areas. Let $\cK_c^n$ denote the space of $K\in\cK_0^n$ that have their centroids at the origin and let $|K|$ denote the $n$-dimensional volume of $K$.
\begin{theorem}\label{aii}
Let $K\in\cK_c^n$ and $B_K\in\cK_c^n$ be the ball such that $|B_K|=|K|$. 
If $\,\phi\in\concave$, then
$$\Omega_{\phi}(K)\le \Omega_{\phi}(B_K)$$
and there is equality for strictly increasing $\phi$  if and only if $K$ is an ellipsoid. 
\end{theorem}
\noindent
For $\phi(t)=t^{1/(n+1)}$ and smooth convex bodies, Theorem \ref{aii} is the classical affine isoperimetric inequality of differential geometry. For general convex bodies, proofs of the classical affine isoperimetric inequality were given by Leichtwei\ss \  \cite{Leichtweiss88}, Lutwak \cite{Lutwak91}, and Hug \cite{Hug96}.  For $L_p$~affine surface areas, the affine isoperimetric inequality was established by Lutwak \cite{Lutwak96} for $p>1$ and by Werner and Ye \cite{Werner:Ye} for $p>0$.

Polarity on convex bodies  induces the following duality on $L_{\phi}$~affine surface areas. Let $K^*=\{x\in\R^n:  x \cdot y \le 1$ for $ y\in K\}$ denote the polar body of $K\in\cK_0^n$. For  $\phi\in\concave$, define $\phi_*:(0,\infty)\to(0,\infty)$ by $\phi_*(s)=s\,\phi(1/s)$.
\begin{theorem}\label{dual}
If $\,\phi\in\concave$, then 
$\,\Omega_{\phi}(K^*) = \Omega_{\phi_*}(K)$ holds
for every $K\in\cK_0^n$.
\end{theorem}
\noindent
For $L_p$~affine surface areas and $p>0$,  Theorem \ref{dual} is due to Hug \cite{Hug96b}: $\Omega_p(K^*)= \Omega_{n^2/p}(K)$  for every $K\in\cK_0^n$.

\smallskip
An alternative definition of $L_p$~affine surface area uses integrals of the curvature function $f(K,\cdot)$ over the unit sphere $\sn$ (see \cite{Lutwak96}). This approach can also be used for $L_{\phi}$~affine surface areas.
\begin{theorem}\label{sphere}
If $\,\phi\in\concave$, then 
$$\,\Omega_{\phi}(K) = \int_{\sn} \phi_*(\a(K,u))\,d\nu_K(u)$$
for every $K\in\cK_0^n$.
\end{theorem}
\noindent
Here $\a(K,u)=\f(K,u)= h(K,u)^{n+1}\,f(K,u)$ is the $L_p$~curvature function of $K$ (see \cite{Lutwak96}) for $p=-n$, while $h(K,u)$ is the support function of $K$, and $d\nu_K(u)= d\cH(u)/h(K,u)^n$ (see Section \ref{tools} for precise definitions). For $L_p$~affine surface areas and $p>0$, Theorem \ref{sphere} is due to Hug \cite{Hug96}.

\smallskip
The family of $L_{\phi}$~affine surface areas for $\phi\in\concave$ includes all $\sln$ invariant and upper semicontinuous valuations on $\cK_0^n$ that vanish on polytopes and, in particular, all $L_p$~affine surface areas for $p>0$. However, $L_p$~affine surface areas for $p<0$ do not be belong to the family of $L_{\phi}$~affine surface areas.
Recent results by Meyer and Werner \cite{Meyer:Werner2000}, Sch\"utt and Werner \cite{Schuett:Werner2004}, 
Werner \cite{Werner2007a}, and Werner and Ye \cite{Werner:Ye} underline the importance of $L_p$~affine surface area also for $p<0$.

A new family of affine surface areas generalizes $L_p$~affine surface area for $-n<p<0$.
Let $\convex$ be the set of functions $\psi:(0,\infty)\to (0,\infty)$ such that $\psi$ is convex, $\lim_{t\to 0}\psi(t)=\infty$, and $\lim_{t\to\infty} \psi(t)=0$. Set $\psi(0)=\infty$. For $\psi\in\convex$, we define the $L_{\psi}$~affine surface area of $K$  by
\begin{equation}\label{two}    
\Omega_{\psi}(K) = \int_{\bd K} \psi (\k(K,x))\,d\mu_K(x).
\end{equation}
The following theorem establishes basic properties of $L_{\psi}$~affine surface areas.

\begin{theorem}\label{ls}
If $\,\psi\in\convex$, then $\Omega_{\psi}(K)$
is positive for every $K\in\cK_0^n$ and  $\Omega_{\psi}(P)=\infty$ for every $P\in\cP_0^n$. In addition,  $\Omega_{\psi}:\cK_0^n\to (0,\infty]$ is both lower semicontinuous and an $\,\sln$ invariant valuation.
\end{theorem}

\noindent
An immediate consequence of Theorem~\ref{ls} is the following result for $L_p$~affine surface area.

\begin{corollary}
If $-n<p<0$, then $\,\Omega_p:\cK_0^n\to(0,\infty]$ 
is positive for every $K\in\cK_0^n$ and  $\Omega_p(P)=\infty$ for every $P\in\cP_0^n$. In addition,  $\,\Omega_p:\cK_0^n\to(0,\infty]$ is  both  lower semicontinuous and an $\,\sln$ invariant valuation.
\end{corollary}

Affine isoperimetric inequalities for $L_{\psi}$~affine surface areas are established in 
\begin{theorem}\label{aii2}
Let $K\in{\cal K}_c^n$ and $B_K\in\cK_c^n$ be the ball such that  $|B_K|=|K|$. 
If $\,\psi\in\convex$, then
$$\Omega_{\psi}(K)\ge \Omega_{\psi}(B_K)$$
and there is equality for strictly decreasing $\psi$ if and only if $K$ is an ellipsoid. 
\end{theorem}
\noindent
For $\psi(t)=t^{p/(n+p)}$ and $-n<p<0$, this result was proved (in a different way) 
by Werner and Ye \cite{Werner:Ye}.

\smallskip
For $\,\psi\in\convex$, define $\Omega^*_{\psi}:\cK_0^n\to(0,\infty]$ by $\Omega^*_{\psi}(K):=\Omega_{\psi}(K^*)$. The following theorem establishes basic properties of these affine surface areas.

\begin{theorem}\label{dual2}
If $\,\psi\in\convex$, then $\Omega^*_{\psi}(K)$
is positive for every $K\in\cK_0^n$ and  $\Omega^*_{\psi}(P)=\infty$ for every $P\in\cP_0^n$. In addition, $\Omega^*_{\psi}:\cK_0^n\to (0,\infty]$ is  both  lower semicontinuous and an $\,\sln$ invariant valuation.
\end{theorem}

\noindent
The family of affine surface areas $\Omega^*_{\psi}$ for $\psi\in\convex$ complements $L_\phi$ affine surface areas for  $\phi\in\concave$ and $L_\psi$ affine surface areas for  $\psi\in\convex$. Whereas $L_\phi$ affine surface areas for  $\phi\in\concave$ include affine surface areas homogeneous of degree $q$ for all $|q|<n$ and $L_\psi$~affine surface areas for  $\psi\in\convex$ include affine surface areas homogeneous of degree $q$ for all $q>n$, the new family includes affine surface areas homogeneous of degree $q$ for all $q<-n$.

The next theorem gives a representation of $\Omega^*_{\psi}$ corresponding to that of Theorem \ref{sphere}.

\begin{theorem}\label{sphere2}
If $\,\psi\in\convex$, then 
$$\,\Omega^*_{\psi}(K) = \int_{\sn} \psi(\a(K,u))\,d\nu_K(u)$$
for every $K\in\cK_0^n$.
\end{theorem}

\smallskip
For $p<-n$, $L_p$~affine surface area was defined by Sch\"utt and Werner \cite{Schuett:Werner2004} using \eqnref{lpdef}. Here a different approach is used and a different definition of $L_p$~affine surface areas for $p<-n$ is given:
\begin{equation}\label{rdef}
\Omega_p(K):=\int_{\sn}\a(K,u)^{\frac{n}{n+p}}\,d\nu_K(u).
\end{equation}
By Theorem \ref{sphere2},  $\Omega_p(K)=\Omega^*_{n^2/p}(K)= \Omega^*_{\psi}(K)$ with $\psi(t)=t^{n/(n+p)}$ and $p<-n$.

An immediate consequence of Theorem~\ref{dual2} is the following result for $L_p$~affine surface area as defined by \eqnref{rdef}.

\begin{corollary}
If $p<-n$, then $\,\Omega_p:\cK_0^n\to(0,\infty]$ 
is positive for every $K\in\cK_0^n$ and  $\Omega_p(P)=\infty$ for every $P\in\cP_0^n$. In addition,  $\,\Omega_p:\cK_0^n\to(0,\infty]$ is  both  lower semicontinuous and an $\,\sln$ invariant valuation.
\end{corollary}

\section{Tools}\label{tools}

Basic notions on convex bodies and their curvature measures are collected. For detailed information, see  \cite{Gardner, Gruber, Schneider:CB}. 
Let $K\in \cK_0^n$.  The support function of $K$ is defined for $x\in\R^n$ by
$$h(K,x)=\max\{ x\cdot y: y\in K\}.$$
The radial function of $K$ is defined for $x\in\R^n$ and $x\neq 0$ by
$$\rho(K,x)=\max\{t>0: t\,x\in K\}.$$
Note that these definitions immediately imply that
\begin{equation}\label{radprop}
\rho(K,x)=1 \fuer x\in\bd K,
\end{equation}
\begin{equation}\label{radhom}
\rho(K, t\,u)= \frac1t \,\rho(K,x)  \fuer t>0,
\end{equation}
and
\begin{equation}\label{polar}
h(K,u)=\frac1{\rho(K^*,u)},
\end{equation}
where $K^*$ is the polar body of $K$.

Let $\cB(\R^n)$ denote the family of Borel sets in $\R^n$ and $\sigma(K,\beta)$ the spherical image of $\beta\in\cB(\R^n)$, that is, the set of all exterior unit normal vectors of $K$ at points of $\beta$. Note that $\sigma(K,\beta)$ is Lebesgue measurable for each $\beta\in\cB(\R^n)$.
For a sequence of convex bodies $K_j\in\cK_0^n$ converging to $K\in\cK_0^n$ and a closed set $\beta\subset\R^n$, we have
\begin{equation}\label{closed}
\limsup_{j\to\infty} \sigma(K_j,\beta)\subset \sigma(K,\beta).
\end{equation}
For $\beta\in\cB(\R^n)$, set 
$$C(K,\beta)=\int_{\sigma(K,\beta)} \frac{d\cH(u)}{h(K,u)^n},$$
where $\cH$ denotes the $(n-1)$-dimensional Hausdorff measure. Hence $C(K,\cdot)$ is a Borel measure on $\R^n$ that is concentrated on $\bd K$.
By \eqnref{polar}, we obtain
\begin{equation}\label{polarvolume}
C(K,\bd K)=n\,|K^*|.
\end{equation}
It follows from \eqnref{closed}  that for every closed set $\beta\subset\R^n$,
\begin{equation}\label{area}
\limsup_{j\to\infty} \,C(K_j,\beta)\le C(K,\beta).
\end{equation}

Let $C_{0}(K,\cdot):\cB(\R^n)\to [0,\infty)$  be the $0$-th {curvature measure} of the convex body $K$ (see \cite{Schneider:CB}, Section 4.2).
For $\beta\in\cB(\R^n)$, we have
\begin{equation}\label{c0H}
C_0(K,\beta)=\cH(\sigma(K,\beta)).
\end{equation}
We decompose the measure $C_0(K,\cdot)$  into measures absolutely continuous and singular with respect to $\cH$, say, $C_0(K,\cdot)=C_0^a(K,\cdot)+C_0^s(K,\cdot)$. Note that
\begin{equation}\label{dev}
\frac{dC_0^a(K,\cdot)}{d\cH}=\kappa(K,\cdot).
\end{equation}

Let $\reg K$ denote the set of regular boundary points of $K$, that is, boundary points with a unique exterior unit normal vector. 
From \eqnref{c0H}, we obtain for $\omega \subset \reg K$ and $\omega\in\cB(\R^n)$,
\begin{equation}\label{reg}
C(K,\omega)=\int_{\sigma(K,\omega)} \frac{d\cH(u)}{h(K,u)^n}=\int_{\omega} \frac{dC_0(K,x)}{(x\cdot u(K,x))^n}.
\end{equation}
We decompose the measure $C(K,\cdot)$  into measures absolutely continuous and singular with respect to the  measure $\mu_K$, say, $C(K,\cdot)=C^a(K,\cdot)+C^s(K,\cdot)$. 
The singular part is concentrated on a $\mu_K$ null set $\omega_0\subset \bd K$, that is, for $\beta\in\cB(\R^n)$
\begin{equation}\label{singular}    
C^s(K,\beta\backslash \omega_0)=0.
\end{equation}
Since $C^a(K,\cdot)$ is concentrated on $\reg K$, \eqnref{dev} and \eqnref{reg} imply 
for $\omega\subset\bd K$ and $\omega\in\cB(\R^n)$,
\begin{equation}\label{absolut}
C^a(K,\omega )= \int_{\omega} \frac{\kappa(K,x)}{(x\cdot u(K,x))^n}\,d\cH(x)=\int_{\omega} \k(K,x)\,d\mu_K(x).
\end{equation}
Combined with \eqnref{polarvolume}, this implies
\begin{equation}\label{dualvolume}
\int_{\partial K} \k(K,x)\,d\mu_K(x)\le n\,|K^*|.
\end{equation}

Hug \cite{Hug96b} proved that for almost all $x\in\bd K$,
$$\kappa(K,x)= \big(\frac{x}{|x|}\cdot u_K(x)\big)^{n+1} f(K^*,\frac{x}{|x|}).$$
Hence we have  for all most all $y\in \bd K^*$,
\begin{equation}\label{trans}
\k(K^*,y)= \a(K, \frac{y}{|y|}).
\end{equation}
Here $|x|$ denotes the length of $x$.

\newpage

\section{Proof of Theorems \ref{aii} and \ref{aii2}} 

Let $\phi\in\concave$ be strictly increasing and $K\in\cK_c^n$. 
By definition \eqnref{one}, Jensen's inequality, \eqnref{dualvolume}, and the monotonicity of $\phi$,  we obtain
\begin{eqnarray*}
\Omega_{\phi}(K)&=& \int_{\partial K} \phi (\k(K,x))\,d\mu_K(x)\\
&\le& n\,|K|\,\, \phi\big( \frac1{n\,|K|}\,\int_{\partial K} \k(K,x)\,d\mu_K(x)\big)\\
&\le& n\,|K|\,\, \phi\big( \frac{|K^*|}{|K|}\big).
\end{eqnarray*}
For origin-centered ellipsoids, $\k(K,\cdot)$ is constant  and there is equality in the above inequalities.
Now we use the Blaschke-Santal\'o inequality: for $K\in\cK_c^n$ 
\begin{equation*}\label{BS}
|K|\,|K^*|\le |B^n|^2
\end{equation*}
with equality precisely for origin-centered ellipsoids (see, for example, \cite{Lutwak1985a}). Here $B^n$ is the unit ball in $\R^n$. 
We obtain
\begin{equation}\label{aiibs}
\Omega_{\phi}(K)\le n\,|K|\,\, \phi\big( \frac{|K^*|}{|K|}\big)\le n\,|K|\,\, \phi\big( \frac{|B^n|^2}{|K|^2}\big)=\Omega_{\phi}(B_K).
\end{equation}
Since $\phi$ is strictly increasing, equality in the second inequality of \eqnref{aiibs} holds if and only if there is equality in the Blaschke-Santal\'o inequality, that is, precisely for ellipsoids. This completes the proof of Theorem \ref{aii} and the proof of Theorem \ref{aii2} follows along similar lines.

\section{Proof of Theorems \ref{dual} and \ref{dual2}}\label{dualp}

Define $\Omega^*_{\phi}$ on $\cK_0^n$ by $\Omega^*_{\phi}(K):= \Omega^{\phantom{*}}_{\phi}(K^*)$. 
Since $\Omega_{\phi}$ is upper semicontinuous, so is $\Omega^*_{\phi}$.
For $K,L, K\cup L\in \cK_0^n$, we have
$$(K\cup L)^* = K^* \cap L^*\,\,\,\hbox{  and }\,\,\, (K\cap L)^* = K^* \cup L^*.$$
Since $\Omega_{\phi}$ is a valuation, this implies that
$$\begin{array}{rcccl}
\Omega^*_{\phi}(K)+\Omega^*_{\phi}(L)&\!\!=
&\!\! \Omega^{\phantom{*}}_{\phi}(K^*)+ \Omega^{\phantom{*}}_{\phi}(L^*)&&\\
&\!\!=&\!\!\Omega^{\phantom{*}}_{\phi}(K^*\cup L^*)+\Omega^{\phantom{*}}_{\phi}(K^*\cap L^*)&&\\
&\!\!=&\!\!\Omega^{\phantom{*}}_{\phi}((K\cap L)^*)+\Omega^{\phantom{*}}_{\phi}((K\cup L)^*)\\
&\!\!=&\!\!\Omega^*_{\phi}(K\cap L)+\Omega^*_{\phi}(K\cup L),
\end{array}$$
that is, $\Omega^*_{\phi}$ is a valuation on $\cK_0^n$.
For  $A\in \sln$ and $K\in\cK_0^n$, we have
$(A\,K)^*= A^{-t}\,K^*$, where $A^{-t}$ denotes the inverse of the transpose of $A$. Since $\Omega_{\phi}$ is $\sln$ invariant,  this implies $\Omega^*_{\phi}(A K)= \Omega^*_{\phi}(K)$, 
that is, $\Omega^*_{\phi}:\cK_0^n\to\R$ is $\sln$ invariant. Since $\Omega_{\phi}$ vanishes on polytopes, so does $\Omega^*_{\phi}$. Therefore $\Omega^*_{\phi}$ satisfies the assumptions of Theorem~\ref{class}. 
Thus there exists $\alpha\in\concave$ such that $\Omega^*_{\phi}=\Omega^{\phantom{*}}_{\alpha}$. Let $\B$ denote the unit ball in $\R^n$. For $r>0$, we obtain from \eqnref{one} that 
\begin{equation*}
\Omega_{\alpha}(r\B)=n\,|\B|\,r^n\,\alpha(\frac1{r^{2n}})
\end{equation*} 
and
\begin{equation*}
\Omega^*_{\phi}(r\B)=\Omega^{\phantom{*}}_{\phi}(\frac1r\B)=\frac{n\,|\B|}{r^n}\, \phi(r^{2n}).  
\end{equation*} 
This shows that $\alpha=\phi_*$ and completes the proof of Theorem \ref{dual}. The proof of Theorem \ref{dual2} follows along the lines of  the proof that $\Omega^*_{\phi}$ satisfies the assumptions of Theorem~\ref{class}.

\section{Proofs of Theorems \ref{sphere} and \ref{sphere2}}

Define $y:\sn\to \bd K^*$ by $u\mapsto \rho(K^*,u)\,u$. Note that this is a Lipschitz function. For the Jacobian $J y$ of $y$, we have  a.e.\ on $\sn$,
\begin{equation}\label{jac}
Jy(u)= \frac{\rho(K^*,u)^{n-1}}{u\cdot u_{K^*}(\rho(K^*,u)\,u)}
\end{equation}
(see, for example, \cite{Hug96b}). By the area formula (see, for example, \cite{Federer}), we have for every a.e.\ defined function $g:\sn \to [0,\infty]$, 
$$
\int_{\sn} g(u) \,Jy(u)\,d\cH(u)= \int_{\bd K^*} g(\frac{y}{|y|})\,d\cH(y).
$$
Setting
$$g(u)= \frac{\tau(\a(K,u))} {h(K,u)^n J y(u)}$$
for $\tau:[0,\infty]\to[0,\infty]$,
we get by   \eqnref{radprop}, \eqnref{radhom}, \eqnref{polar}, and  \eqnref{trans},
\begin{eqnarray*}
\int\limits_{\sn} \!\!\tau(\a(K,u))\,d\nu_K(u)\!\!\!&=&\!\!\!\!\int\limits_{\sn} \!\tau(\a(K,u))\,\frac{d\cH(u)}{h(K,u)^n}\\
&=&\!\!\!\!\int\limits_{\bd K^*}\! \tau(\k(K^*,y))\,\frac{\frac{y}{|y|}\cdot u_{K^*}(y)}{\rho(K^*,\frac{y}{|y|})^{n-1}}\,\rho(K^*,\frac{y}{|y|})^n\,d\cH(y)\\
&=&\!\!\!\!\int\limits_{\bd K^*}\! \tau(\k(K^*,y))\,d\mu_{K^*}(y).
\end{eqnarray*}
For $\tau\in\convex$, this implies Theorem \ref{sphere2}. To obtain Theorem \ref{sphere}, we set $\tau=\phi_*\in\concave$ and apply Theorem \ref{dual}.

\section{Proof of Theorem \ref{ls}}

Let $\psi\in\convex$ and $K\in\cK_0^n$. Note that $\psi$ is strictly decreasing and positive. By definition \eqnref{two}, the Jensen inequality, \eqnref{dualvolume}, and the monotonicity of $\psi$,  we obtain
\begin{eqnarray*}
\Omega_{\psi}(K)&=&\int_{\bd K} \psi(\k(K,x))\,\,d\mu_K(x)\\
 &\ge& n\,|K|\,\, \psi\big( \frac1{n\,|K|}\,\int_{\partial K} \k(K,x)\,d\mu_K(x)\big)\\
&\ge& n\,|K|\,\, \psi\big( \frac{|K^*|}{|K|}\big).
\end{eqnarray*}
This shows that $\Omega_{\psi}(K)>0$. The $\sln$ invariance of $\Omega_{\psi}$ follows immediately from the definition. So does the fact that $\Omega_{\psi}(P)=\infty$ for $P\in\cP_0^n$.

Next, we show that $\Omega_{\psi}$ is a valuation on $\cK_0^n$, that is, for $K,L\in\cK_0^n$  such that $K\cup L\in\cK_0^n$,
\begin{equation}\label{valuation}
\as(K\cup L)+\as(K\cap L)=\as(K)+\as(L).
\end{equation}
Let $K^c=\{x\in\R^n: x\not\in K\}$ and let $\inte K$ denote the interior of $K$. We follow Sch\"utt \cite{Schuett93} (see also \cite{Haberl:star})
and work with the decompositions 
\begin{eqnarray*}
\bd (K\cup L)&=&\big(\bd K\cap \bd L\big) \cup  \big(\bd K\cap L^c\big)\cup \big(\bd  L\cap K^c\big),\\
\bd (K\cap L)&=&\big(\bd K\cap \bd L\big)  \cup \big(\bd K\cap \inte L\big) \cup \big(\bd  L\cap \inte K\big),\\
\bd K&=&\big(\bd K\cap \bd L\big) \cup \big(\bd  K\cap L^c\big) \cup \big(\bd K\cap \inte L\big),\\
\bd L&=&\big(\bd K\cap \bd L\big) \cup \big(\bd  L\cap K^c\big) \cup \big(\bd L\cap \inte K\big),
\end{eqnarray*}
where all unions on the right hand side are disjoint.  Note that for $x$ such that the curvatures $\k(K,x)$, $\k(L,x)$, $\k(K\cup L,x)$, and 
$\k(K\cap L,x)$ exist, 
\begin{equation}\label{cone}
u(K,x)=u(L,x)=u(K\cup L,x)=u(K\cap L,x)
\end{equation}
and
\begin{equation}\label{curv}
\begin{array}{rcl}
\k(K\cup L,x)&=&\min\{\k(K,x), \k(L,x)\},\\
\k(K\cap L,x)&=&\max\{\k(K,x), \k(L,x)\}.\\
\end{array}
\end{equation}
To prove \eqnref{valuation}, we use \eqnref{two}, split the involved integrals using the above decompositions, and use  \eqnref{cone} and \eqnref{curv}.

Finally, we show that $\as$ is lower semicontinuous on $\cK_0^n$. The proof  complements the proofs in \cite{Ludwig:curvature} and  \cite{Lutwak92}.
Let $K\in\cK_0^n$  and $\varepsilon>0$ be chosen.
Since $\k(K,\cdot)$ is measurable a.e.~on $\bd K$ and since the set  $\omega_0$, where the singular part of $C(K,\cdot)$ is concentrated, is a $\mu_K$ null set, we can choose by Lusin's theorem (see, for example, \cite{Federer}) 
pairwise disjoint closed sets $\omega_l\subset \partial K$, $l=1,\ldots$, such that $\k(K,\cdot)$ is continuous as a function restricted to $\omega_l$, such that for every $l=1,\ldots$,
\begin{equation}\label{null}    
\omega_l\cap \omega_0=\emptyset
\end{equation} and such that
\begin{equation}\label{lusin}
\mu_K(\bigcup_{l=1}^\infty \omega_l)=\mu_K(\bd K).
\end{equation}
For $\omega\subset \R^n$, let $\bar \omega$ be the cone generated by $\omega$, that is, $\bar \omega =\{t\,x\in\R^n: t\ge0, x\in \omega\}$. Note that $\bar\omega_l$ is closed and that $\bd K\cap \bar \omega_l=\omega_l$.
 
Let $K_j$ be a sequence of convex bodies converging to $K$. First, we show that for $l=1,\ldots$,
\begin{equation}\label{main}
\liminf_{j\to\infty} \int_{\bd K_j\cap \bar\omega_l} \psi(\k(K_j,x))\,d\mu_{K_j}(x)\ge \int_{\bd K\cap \bar\omega_l} \psi(\k(K,x)\,d\mu_K(x).
\end{equation}
Let $\eta>0$ be chosen. We choose a monotone sequence  $t_i\in(0,\infty)$,  $i=\Z$, $\lim_{i\to -\infty} t_i=0$, $\lim_{i\to \infty} t_i=\infty$,
such that 
\begin{equation}\label{fein}    
\max_{i\in\Z}|\psi(t_{i+1})-\psi(t_i)|\le \eta
\end{equation}
and such that for $i\in\Z$, $j\ge 0$, 
\begin{equation}\label{gut}
\mu_{K_j}(\{x\in \bd K_j : \k(K_j,x)=t_i\})=0,
\end{equation}
where $K_0=K$. This is possible, since $\mu_{K_j}(\{x\in K_j: \k(K_j,x)=t\})>0$ holds only for countably many $t$. Set
$$\omega_{li}=\{x\in\omega_l: t_{i}\le \k(K,x)\le t_{i+1}\}.$$
Since $\k(K,\cdot)$ is continuous on $\omega_l$ and $\omega_l$ is closed, the sets  $\bar \omega_{li}$  are  closed for $i\in\Z$. This implies  by (\ref{area}) that
\begin{equation}\label{Cconv}
\limsup_{j\to\infty} C(K_j,\bar\omega_{li})\le C(K,\bar\omega_{li}).
\end{equation}
By (\ref{null}),  (\ref{singular}), and the definition of $\omega_{li}$, 
\begin{equation}\label{dom}
C(K,\bar\omega_{li})=C^a(K,\bar\omega_{li}) \le t_{i+1}\,\mu_K(\bd K\cap\bar\omega_{li}).
\end{equation}
By (\ref{absolut}),
\begin{equation}\label{abs}
\int_{\bd K_j\cap \bar\omega_{li}}\k(K_j,x)\,d\mu_{K_j}(x)\le C(K_j,\bar\omega_{li}).
\end{equation}
Using the monotonicity of $\psi$, we  obtain
\begin{equation}\label{erstens}
\begin{array}{rcl}
\displaystyle\int_{\omega_l} \psi(\k(K,x))\,d\mu_K(x)
&\le&\displaystyle\sum_{i\in\Z}\,\int_{\omega_{li}}\psi(\k(K,x))\,d\mu_K(x)\\
\displaystyle&\le&\displaystyle\sum_{i\in\Z}\,\psi(t_{i})\,\mu_K(\omega_{li}).
\end{array}
\end{equation}
Using \eqnref{gut}, the Jensen inequality, \eqnref{abs}, and the monotonicity of $\psi$, we obtain
\begin{eqnarray*}
\int\limits_{\bd K_j\cap\bar\omega_l} \!\psi(\k(K_j,x))\, d\mu_{K_j}(x)
&=&\sum_{i\in\Z} \,\,\int\limits_{\bd K_j\cap \bar\omega_{li}} \psi(\k(K_j,x))\, d\mu_{K_j}(x)\\
&=&\sum_{i\in\Z}\!\!{\phantom{|}^\prime}\int\limits_{\bd K_j\cap \bar\omega_{li}} \psi(\k(K_j,x))\, d\mu_{K_j}(x)\\
&\ge&\sum_{i\in\Z}\!\!{\phantom{|}^\prime} \psi
\left(\frac{C(K_j,\bar\omega_{li})}{\mu_{K_j}(\bd K_j\cap\bar\omega_{li})}\right)\,\mu_{K_j}(\bd K_j\cap\bar\omega_{li})
\end{eqnarray*}
where the $'$ indicates that we sum only over $\bar\omega_{li}$ with $\mu_{K_j}(\bd K_j\cap\bar\omega_{li})\ne 0$.
Since 
\begin{eqnarray*}
\liminf_{j\to\infty} &&\hspace{-0.75cm} \sum_{i\in\Z}\!\!{\phantom{|}^\prime} \psi
\left(\frac{C(K_j,\bar\omega_{li})}{\mu_{K_j}(\bd K_j\cap\bar\omega_{li})}\right)\,\mu_{K_j}(\bd K_j\cap\bar\omega_{li})\\
&\ge&
\sum_{i\in\Z}\!\!{\phantom{|}^\prime}  \psi
\left(\limsup_{j\to\infty} \left(\frac{C(K_j,\bar\omega_{li})}{\mu_{K_j}(\bd K_j\cap\bar\omega_{li})}\right)\right)\,\liminf_{j\to\infty} \mu_{K_j}(\bd K_j\cap\bar\omega_{li}),
\end{eqnarray*}
we obtain by \eqnref{Cconv}, \eqnref{dom},(\ref{erstens}),(\ref{fein}), and \eqnref{gut} that
\begin{eqnarray*}
\liminf_{j\to\infty} &&\hspace{-0.95cm}\int_{\bd K_j\cap \bar\omega_{l}} \psi(\k(K_j,x))\,d\mu_{K_j}(x)\\
&\ge&\sum_{i\in\Z}\!\!{\phantom{|}^\prime} 
\psi\left( \frac{C(K,\bar\omega_{li})}{\mu_K(\bd K\cap\bar\omega_{li})}\right)\,\mu_K(\bd K\cap\bar\omega_{li})\\
&\ge&\sum_{i\in\Z}\,\psi (t_{i+1})\,\mu_K(\bd K\cap\bar\omega_{li})\\
&=&\sum_{i\in\Z}\,\psi (t_{i})\,\mu_K(\bd K\cap\bar\omega_{li})-\sum_{i\in\Z}\left(\psi (t_{i})-\psi (t_{i+1})\right)\,\mu_K(\bd K\cap\bar\omega_{li})\\
&\ge&\int_{\bd K\cap \bar\omega_l} \psi (\k (K,x))\,d\mu_K(x)- \eta \,\mu_K(\bd K \cap \bar \omega_l).
\end{eqnarray*}
Since $\eta>0$ is arbitrary, this proves (\ref{main}).

Finally, \eqnref{gut} and \eqnref{main} imply
\begin{eqnarray*}
\liminf_{j\to\infty} \int\limits_{\bd K_j} \psi (\k (K_j, x))\,d\mu_{K_j}(x)
&=& \liminf_{j\to\infty} \sum_{l=1}^\infty \int\limits_{\bd K_j\cap\bar \omega_l} \psi (\k (K_j, x))\,d\mu_{K_j}(x)\\
&\ge& \sum_{l=1}^\infty \liminf_{j\to\infty} \,\int\limits_{\bd K_j\cap\bar \omega_l}\! \psi ( \k (K_j, x))\,d\mu_{K_j}(x)\\
& \ge & \int\limits_{\bd K} \psi (\k (K, x))\,d\mu_{K}(x).
\end{eqnarray*}
This completes the proof of the theorem.

\section{Open problems}\label{conj}

The affine surface areas $\,\Omega_{\psi}$ and $\Omega^*_{\psi}$ for $\psi\in\convex$ are lower semicontinuous and $\sln$ invariant valuations. 
More general examples of such functionals are
$$\Psi= \Omega^{\phantom{*}}_{\psi_1} + \Omega^{*}_{\psi_2} -\Omega^{\phantom{*}}_{\phi}$$
for $\psi_1, \psi_2\in\convex$ and $\phi\in \concave$. Additional examples are the following continuous functionals
\begin{equation*}\label{cont}
|K|\mapsto c_0 + c_1\,|K|+c_2\,|K^*|
\end{equation*}
for $c_0,c_1,c_2\in\R$.
In view of  Theorem~\ref{class}, this gives raise to the following

\begin{question}
If $\,\Psi:\cK_0^n\to (-\infty,\infty]$ is a lower semicontinuous and $\,\sln$ invariant valuation, then there exist $\psi_1, \psi_2\in\convex$, $\phi\in\concave$, and $c_0,c_1,c_2\in\R$ such that 
$$\Psi(K)= c_0 + c_1\,|K|+c_2\,|K^*| +\Omega^{\phantom{*}}_{\psi_1}(K) + \Omega^{*}_{\psi_2}(K) -\Omega^{\phantom{*}}_{\phi}(K)$$
for every $K\in\cK_0^n$.
\end{question}

The following special case of the above conjecture is of particular interest.

\begin{question}
If $\,\Psi:\cK_0^n\to (-\infty,\infty]$ is a lower semicontinuous and $\,\sln$~invariant valuation  that is homogeneous of degree $q<-n\,$ or $q>n$, then there exists $c\ge 0$ such that
$$\Psi(K)= c\, \Omega_p(K)$$
for every $K\in\cK_0^n$, where $p=n\,(n-q)/(n+q)$.
\end{question}

\parindent=0pt

\normalsize
\bigskip
\begin{samepage}
Department of Mathematics, Polytechnic Institute of New York University, 6 MetroTech Center, Brooklyn, NY 11201, U.S.A.\\

E-mail: mludwig@poly.edu
\end{samepage}

\end{document}